\documentclass[12pt]{amsart}
\usepackage{amsmath,amssymb,amsfonts,amsthm,latexsym}

\DeclareMathAlphabet{\curly}{U}{rsfs}{m}{n}

\newtheorem{lem}{Lemma}[section]

\newcommand{\ZZ}{{\mathbb Z}}

\newcommand{\RR}{{\mathbb R}}

\newcommand{\NN}{{\mathbb N}}
\newcommand{\bb}{{\mathbf b}}

\newcommand{\LL}{\curly L}

\newcommand{\BB}{\curly B}

\newcommand{\TT}{\curly T}
\newcommand{\PP}{\curly P}

\newcommand{\Vol}{\operatorname{Vol}}   

\newcommand{\g}{\ensuremath{\gamma}}
\newcommand{\del}{\ensuremath{\delta}}
\newcommand{\lam}{\ensuremath{\lambda}}

\newcommand{\zz}{\ensuremath{\mathbf{z}}}

\newcommand{\bj}{\ensuremath{\mathbf{j}}}

\newcommand{\bx}{{\ensuremath{\boldsymbol{\xi}}}}

\newcommand{\fl}[1]{{\ensuremath{\left\lfloor {#1} \right\rfloor}}}

\newcommand{\pfrac}[2]{{\left(\frac{#1}{#2}\right)}}

\newcommand{\be}{\begin{equation}}
\newcommand{\ee}{\end{equation}}
\newcommand{\benn}{\begin{equation*}}   
\newcommand{\eenn}{\end{equation*}}

\newcommand{\EE}{\mathbb E}

\renewcommand{\AA}{\curly A}
\renewcommand{\(}{\left(}
\renewcommand{\)}{\right)}

\numberwithin{equation}{section}

\newif\ifdraft

\begin{document}

\title{Integers with a divisor in $(y,2y]$}

\author{Kevin Ford}

\address{Department of Mathematics, 1409 West Green Street, University
of Illinois at Urbana-Champaign, Urbana, IL 61801, USA}
\email{ford@math.uiuc.edu}

\date{\today}
\thanks{Research supported by National Science Foundation grants
DMS-0301083 and DMS-0555367.}

\begin{abstract} We determine, up to multiplicative constants, how many
  integers $n\le x$ have a divisor in $(y,2y]$.
\end{abstract}

\maketitle

%
%
%
\section{Introduction}\label{sec:intro}
%
%
%
Let $H(x,y,z)$ be the number of integers $n\le x$ which have a divisor in
the interval $(y,z]$.  In the author's paper \cite{F}, the correct order of
growth of $H(x,y,z)$ was determined for all $x,y,z$.  In particular,
\be\label{Hxy2y}
H(x,y,2y) \asymp \frac{x}{(\log y)^{\delta} (\log\log y)^{3/2}} \qquad
(3\le y\le \sqrt{x}),
\ee
where
$$
\del = 1 - \frac{1+\log\log 2}{\log 2} = 0.086071\ldots.
$$

In this note we prove only the important special case \eqref{Hxy2y}, 
omitting the parts of the argument required for other cases. 
In addition, we present an alternate proof, dating from 2002, of the 
lower bound implicit in \eqref{Hxy2y}.
This proof avoids the use of results about uniform order statistics required in
\cite{F}, and instead utilizes the cycle lemma from combinatorics.
Although shorter and 
technically simpler than the argument in \cite{F}, this method
is not useful for a related problem also considered in \cite{F}, that of
counting integers with a prescribed number of divisors in $(y,2y]$.
We also simplify the upper bound argument using  a result on
sums of arithmetic functions due to Koukoulopoulos \cite[Lemma
  2.2]{DK}, a short proof of which we give below.

We mention here one of the applications of \eqref{Hxy2y}, a 1955 problem
of Erd\H os (\cite{Erdos55}, \cite{Erdos60}) 
known colloquially as the ``multiplication table problem''.
Let $A(x)$ be the number of positive integers $n\le x$ which can be written
as $n=m_1 m_2$ with each $m_i \le \sqrt{x}$.   Then
$$
A(x) \asymp \frac{x}{(\log x)^{\del} (\log\log x)^{3/2}}.
$$
This follows directly from \eqref{Hxy2y} and the inequalities
$$
H\( \frac{x}{4}, \frac{\sqrt{x}}{4}, \frac{\sqrt{x}}{2} \) \le A(x) \le
\sum_{k\ge 0} H\(\frac{x}{2^k}, \frac{\sqrt{x}}{2^{k+1}}, 
\frac{\sqrt{x}}{2^k} \). 
$$

More on the history of estimations of $H(x,y,z)$, further applications
and references may be found in \cite{F}.

\bigskip

{\bf Heuristic argument.}
For brevity, let $\tau(n,y,z)$ be the number of divisors of $n$ in $(y,z]$.
Write $n=n'n''$, where $n'$ is composed only of primes $\le 2y$ and
$n''$ is composed only of primes $>2y$.
For simplicity, assume $n'$ is squarefree and $n'\le y^{100}$.
Assume for the moment that the set $D(n')=\{ \log d : d|n' \}$ is
uniformly distributed in $[0,\log n']$.  If $n'$ has $k$ prime factors, then 
the expected value of $\tau(n',y,2y)$ should be about 
$\frac{2^k\log 2}{\log n'} \asymp
\frac{2^k}{\log y}$.  This is $\gg 1$ precisely when $k\ge k_0 + O(1)$,
where $k_0 := \fl{\frac{\log\log y}{\log 2}}$.
Using the fact (e.g. Theorem 08 of \cite{Divisors}) 
that the number of $n\le x$ with $n'$ having $k$ prime factors is
of order
$$
\frac{x}{\log y} \frac{(\log\log y)^k}{k!},
$$
we obtain a heuristic estimate for $H(x,y,2y)$ of order
$$
\frac{x}{\log y} \sum_{k\ge k_0+O(1)} \frac{(\log\log y)^k}{k!} 
\asymp \frac{x (\log\log y)^{k_0}}{k_0! \log y} 
\asymp \frac{x}{(\log y)^{\del} (\log\log y)^{1/2}}.
$$

This is slightly too big, and the reason stems from the uniformity
 assumption about $D(n')$.  In fact, for most $n'$ with about $k_0$ prime
 factors, the set $D(n')$ is far from uniform, possessing many clusters
of divisors and large gaps between clusters.  This substantially decreases
the likelihood that $\tau(n',y,2y)\ge 1$.  The numbers $\log\log p$ over $p|n'$
are well-known to behave like random numbers in $[0, \log\log 2y]$.
Consequently, if we write $n'=p_1\cdots p_k$, where $p_1<p_2<\ldots<p_k$, then 
we expect $\log\log p_j \approx \frac{j\log\log y}{k_0} = j\log 2+O(1)$ for
each $j$.  Large deviation results from probability theory
(see Smirnov's theorem in \S \ref{sec:upper2}; also see Ch. 1 of
\cite{Divisors}) tell us that with high probability
there is a $j$ for which $\log\log p_j \le j\log 2 - c\sqrt{\log\log y}$,
where $c$ is a small positive constant.  Thus, the $2^j$ divisors of
$p_1\cdots p_j$ will be clustered in an interval of logarithmic length about
$\ll \log p_j \le 2^j e^{-c\sqrt{\log\log y}}$.  On a logarithmic scale, the
divisors of $n'$ will then lie in $2^{k-j}$ translates of this cluster.
A measure of the degree of clustering of the divisors of an integer 
$a$ is given by
$$
L(a)=\text{meas} \LL(a), \qquad \LL(a)=\bigcup_{d|a} [-\log 2+\log d,\log d).
$$
The probability that $\tau(n',y,2y)\ge 1$ should then be about $L(n')/\log y$.
Making this precise leads to the upper and lower bounds for $H(x,y,2y)$
given below in Lemmas \ref{HL_lower} and \ref{HL}.  The upper bound for $L(a)$
given in Lemma \ref{L} (iii) below quantifies 
how small $L(a)$ must be when there is 
a $j$ with $\log\log p_j$ considerably smaller than $j\log 2$.

What we really need to count is $n$ for which $n'$ has about $k_0$ 
prime factors \emph{and} $L(n') \gg \log n'$.  This roughly corresponds to 
asking for $\log\log p_j \ge j\log 2 - O(1)$ for all $j$.
The anologous problem
from statistics theory is to ask for the likelihood than given $k_0$ random
numbers in $[0,1]$, there are $\le k_0x+O(1)$ of them which are $\le x$,
uniformly in $0\le x\le 1$.  In section \ref{sec:upper2}, Lemma \ref{Q},
we will see that this probability is about
$1/k_0 \asymp 1/\log\log y$ and this leads to the correct order \eqref{Hxy2y}.

\medskip
\noindent
{\bf Notation:} 
Let $\tau(n)$ be the number of positive divisors of $n$, 
and  define $\omega(n)$ to be the number of distinct prime
divisors of $n$.  Let
$P^+(n)$ be the largest prime factor of $n$ and let $P^-(n)$ be the 
smallest prime factor of $n$.  Adopt the notational conventions $P^+(1)=0$ and
$P^-(1)=\infty$.  Constants implied by $O$, $\ll$ and $\asymp$
are absolute. The notation $f \asymp g$ means $f\ll g$ and $g\ll f$.

We shall make frequent use of the following estimate,
which is a consequence of the Prime Number Theorem with
classical de la Val\'ee Poussin error term.
For certain constants $c_0, c_1$, 
\be\label{sum1p}
\sum_{p\le x} \frac{1}{p} = \log\log x + c_0 + O(e^{-c_1\sqrt{\log x}})
\qquad (x\ge 2).
\ee
We also 
need the standard sieve bound (e.g. \cite{HR}; Theorem 06 and Exercise 02
of \cite{Divisors})
\be\label{sieve}
| \{ n\le x : P^-(n)>z \}| \asymp \frac{x}{\log z} \qquad (x \ge 2z \ge 4)
\ee
and Stirling's formula $k! \sim \sqrt{2\pi k} (k/e)^k$.

%
%
%
%
%
%
\section{Lower bound}\label{sec:lower}
%
%
%
%
%

In this section we prove the lower bound implicit in \eqref{Hxy2y}.
The first step is to bound $H(x,y,2y)$ in terms of a sum of $L(a)/a$.
Next, sums of $L(a)/a$ are related via the Cauchy-Schwarz inequality to
sums of a function $W(a)$ which counts pairs of divisors of 
$a$ which are close together.
With a strategic choice of sets of $a$ to average over, the problem is
reduced to the estimation of a certain combinatorial sum.  This is
accomplished with the aid of a tool closely related to the 
so-called ``cycle lemma''.

\begin{lem}\label{HL_lower}
If $3\le y\le \sqrt{x}$, then
$$
H(x,y,2y) \gg \frac{x}{\log^2 y}\sum_{a\le y^{1/8}} \frac{L(a)}{a}.
$$
\end{lem}

\begin{proof}
Let $y_0$ be a sufficiently large constant.  If $3\le y\le y_0$, then
$H(x,y,2y) \gg x \gg \frac{x L(1)}{\log^2 y}$ since $L(1)=\log 2$.
If $y\ge y_0$, consider integers $n=apb\le x$ with
$a\le y^{1/8}$, all prime factors of $b$ are $>2y$ or in $[y^{1/4},y^{3/4}]$,
and $p$ is a prime with $\log(y/p)\in \LL(a)$.  The last condition implies
that $\tau(ap,y,2y)\ge 1$.  In particular,
$y^{7/8} \le y/a < p \le 2y$.  Thus, each $n$ has a unique
representation in this form.  Fix $a$ and $p$ and note that
$x/(ap) \ge x/(2y^{9/8}) \ge \frac12 y^{7/8}$.
  If $x/(ap) \ge 4y$,
\eqref{sieve} implies that the number of $b\le \frac{x}{ap}$ with
$P^-(b)>2y$ is $\gg \frac{x}{ap\log y}$.
If $x/(ap) < 4y$, then the number of
$b\le \frac{x}{ap}$ composed of two prime factors in $(y^{1/4},y^{3/4}]$
is likewise $\gg \frac{x}{ap\log y}$.   Hence
$$
H(x,y,2y) \gg \frac{x}{\log y} \sum_{a\le y^{1/8}} \frac{1}{a} 
\sum_{\log(y/p)\in\LL(a)} \frac{1}{p}.
$$
Since $\LL(a)$
is the disjoint union of intervals of length $\ge \log 2$ and
$p\ge y^{7/8}$, for
each $a$ we have by repeated application of \eqref{sum1p}
\[
\sum_{\log(y/p)\in\LL(a)} \frac{1}{p} \gg \frac{L(a)}{\log y}. \qedhere
\]
\end{proof}

\begin{lem}\label{LW}
For any finite set $\AA$ of positive integers,
$$
\sum_{a\in \AA} \frac{L(a)}{a} \ge \frac{\( \sum_{a\in \AA} \frac{\tau(a)}{a}
  \)^2}{6 \sum_{a\in \AA} \frac{W(a)}{a} },
$$
where
$$
W(a) = | \{ (d,d') : d|a, d'|a, |\log d/d'| \le \log 2 \}|.
$$
\end{lem}

\begin{proof} Since $\tau(a)\log 2 = \int \tau(a,e^u,2e^u)\, du$, 
by the Cauchy-Schwarz inequality,
\begin{align*}
\( \sum_{a\in\AA} \frac{\tau(a)}{a} \)^2 (\log 2)^2 
&= \( \sum_{a\in\AA} \frac{1}{a} \int \tau(a,e^u,2e^u)\, du \)^2 \\
&\le \( \sum_{a\in \AA} \frac{L(a)}{a} \) \( \sum_{a\in\AA} \frac{1}{a}
\int \tau^2(a,e^u,2e^{u})\, du \).
\end{align*}
Let $k_j=\tau(a,2^{j-1},2^j)$ for each integer $j$.  Then
\[
\int \tau^2(a,e^u,2e^{u})\, du \le (\log 2) \sum_{j} (k_j+k_{j+1})^2 \le
4(\log 2) \sum_j k_j^2 \le 4 (\log 2) W(a). \qedhere
\]
\end{proof}

We apply Lemma \ref{LW} with sets $\AA$ of integers whose prime
factors are localized.  To simplify later analysis, 
partition the primes into  sets $D_1, D_2, \ldots$, where
each $D_j$ consists of the primes in an interval
$(\lam_{j-1},\lam_j]$, with $\lam_j \approx \lam_{j-1}^2$.  More precisely, 
let $\lam_0=1.9$ and define inductively $\lam_j$ for $j\ge 1$ as the
largest prime so that
\be\label{Dj}
\sum_{\lam_{j-1} < p \le \lam_j} \frac{1}{p} \le \log 2.
\ee
For example, $\lam_1=2$ and $\lam_2=7$. 
By \eqref{sum1p}, we have
$$
\log\log \lam_{j} -\log\log \lam_{j-1} = \log 2 + O(e^{-c_1\sqrt{\log
    \lam_{j-1}}}), 
$$
and thus for some absolute constant $K$,
\be\label{lam}
2^{j-K} \le \log \lam_j \le  2^{j+K} \qquad (j\ge 0).
\ee
For a vector $\bb=(b_1,\ldots,b_J)$ of non-negative integers,
let $\AA(\bb)$ be the set of square-free integers $a$ composed of exactly 
$b_j$ prime factors from $D_j$ for each $j$. 

%
%

\begin{lem}\label{sumW}
Assume  $\bb=(b_1,\ldots,b_J)$. Then
$$
\sum_{a\in \AA(\bb)} \frac{W(a)}{a} \ll \frac{(2\log 2)^{b_1+\cdots+b_J}}{b_1! 
  \cdots b_J!}  \sum_{j=1}^J 2^{-j+b_1 + \cdots + b_j}.
$$
\end{lem}

\begin{proof} Let $B=b_1+\cdots+b_J$ and for $j\ge 0$ let $B_j=\sum_{i\le j}
b_j$.  Let $a=p_1\cdots p^{\phantom{p}}_B$, where
\be\label{sumW_a}
p^{\phantom{p}}_{B_{j-1}+1}, \ldots, p^{\phantom{p}}_{B_j} \in D_j \qquad (1\le j\le J)
\ee
and the primes in each interval $D_j$ are unordered. 
Since $W(p_1\cdots p_B)$
is the number of pairs $Y,Z\subseteq \{1,\ldots,B\}$ with
\be\label{sumW_b}
\left| \sum_{i\in Y} \log p_i - \sum_{i\in Z} \log p_i \right| \le \log 2,
\ee
we have
\be\label{sumW_c}
\sum_{a\in \AA(\bb)} \frac{W(a)}{a} \le \frac{1}{b_1! \cdots b_J!}
\sum_{Y,Z \subseteq \{1,\ldots,B\}} \;\; \sum_{\substack{p_1,\ldots,p_B \\
\eqref{sumW_a}, \eqref{sumW_b}}} \frac{1}{p_1 \cdots p_B}.
\ee
When $Y=Z$, \eqref{Dj} implies that 
the inner sum on the right side 
of \eqref{sumW_c} is $\le (\log 2)^B$, and there
are $2^B$ such pairs $Y,Z$.  When $Y\ne Z$, let $I=\max [(Y\cup Z)-(Y\cap
Z)]$.   With all the $p_i$ fixed except for $p_I$, \eqref{sumW_b} implies that
$U \le p_I \le 4U$ for some number $U$.  Let
$E(I)$ be defined by $B_{E(I)-1} <  I \le B_{E(I)}$, i.e. 
$p_I \in D_{E(I)}$.   By  \eqref{sum1p},
$$
\sum_{\substack{U\le p_I \le 4U \\ p_I\in D_{E(I)}}} \frac{1}{p_I} 
\ll \frac{1}{\max(\log U,\log \lam_{E(I)-1})} \ll 2^{-E(I)}.
$$
Thus, by \eqref{Dj} 
the inner sum in \eqref{sumW_c} is $\ll 2^{-E(I)} (\log 2)^B$.  With $I$
fixed, there correspond $2^{B-I+1} 4^{I-1}=2^{B+I-1}$ pairs $Y,Z$.  By
\eqref{sumW_c}, 
$$
\sum_{a\in \AA(\bb)} \frac{W(a)}{a} \ll \frac{(2\log 2)^B}{b_1! \cdots b_J!}
\left[ 1 + \sum_{I=1}^B 2^{I-E(I)} \right]
\ll\frac{(2\log 2)^B}{b_1! \cdots b_J!}
\sum_{j=1}^J 2^{-j} \sum_{B_{j-1}<I\le B_j} 2^I,
$$
and the claimed bound follows.
\end{proof}

Now suppose that $M$ is a sufficiently large positive integer,
$b_i=0$ for $i<M$, and $b_j\le Mj$ for each $j$.  By \eqref{lam},
\be\label{sumtau}
\begin{split}
\sum_{a\in\AA(\bb)} \frac{\tau(a)}{a} &= 2^k \prod_{j=M}^J 
\frac{1}{b_j!} \biggl(
  \sum_{p_1\in D_j} \frac{1}{p_1} \sum_{\substack{p_2\in D_j \\ p_2 \ne p_1}}
  \frac{1}{p_2} \cdots \sum_{\substack{p_{b_j}\in D_j \\ p_{b_j} \not\in
  \{ p_1, \ldots, p_{b_j-1} \} }} \frac{1}{p_{b_j}} \biggr) \\
&\ge 2^k \prod_{j=M}^J \frac{1}{b_j!}
   \biggl( \log 2 - \frac{b_j}{\lam_{j-1}} \biggr)^{b_j} \\
&\ge \frac{(2\log 2)^k}{2 b_M! \cdots b_J!}.
\end{split}
\ee

Let
$$
k = \fl{\frac{\log\log y}{\log 2} - 2M}, \qquad J=M+k-1.
$$
Let $\BB$ be the set of vectors $(b_1,\ldots,b_J)$ with $b_i=0$ for
$i<M$ and $b_1+\cdots+b_J=k$.  Let $\BB^*$ be the set of $\bb\in\BB$
with $b_j\le \min(Mj,M(J-j+1))$ for each $j\ge M$.  If $\bb\in \BB^*$
and $a\in \AA(\bb)$, then by \eqref{lam},
$$
\log a \le \sum_{j=M}^J b_j \log \lam_j \le 
M \sum_{l=0}^{J-M} (l+1) 2^{J+K-l} < \frac{\log y}{8}
$$
if $M$ is large enough, as $2^{J+1} \le 2^{-M}\log y$.  Put
\be\label{fdef}
f(\bb)=\sum_{h=M}^J 2^{M-1-h + b_M+\cdots+b_h}.
\ee
By Lemma \ref{sumW},
$$
\sum_{a\in \AA(\bb)} \frac{W(a)}{a} \ll \frac{(2\log 2)^k}{b_M!\cdots b_J!} \(
1 + 2^{1-M} f(\bb) \) \ll \frac{(2\log 2)^k}{b_M!\cdots b_J!} f(\bb)
$$
since $f(\bb)\ge 1/2$.
By Lemmas \ref{HL_lower} and \ref{LW}, plus \eqref{sumtau}, we have
for large $y$
\be\label{low1}
H(x,y,2y) \gg \frac{x (2\log 2)^k}{\log^2 y} \sum_{\bb\in \BB^*}
\frac{1}{b_M!\cdots b_J! f(\bb)}.
\ee
Observe that the product of factorials is unchanged under permutation
of $b_M,\ldots,b_J$.  Roughly speaking, 
$$
f(\bb) \approx g(\bb) := \max_{j} 2^{(b_M-1)+\cdots+(b_j-1)}.
$$
Note that $(b_M-1)+\cdots+(b_J-1)=k-(J-M+1)=0$.

Given real numbers $z_1,\cdots,z_k$ with zero sum, there is a cyclic
permutation $\zz'$ of the vector $\zz=(z_1,\ldots,z_k)$ all of whose partial
sums are $\ge 0$: let $i$ be the index minimizing $z_1+\cdots+z_i$ and take
$\zz'=(z_{i+1},\ldots,z_k,z_1,\ldots,z_{i})$.  In combinatorics, this fact is
know as the \emph{cycle lemma}.  Thus, there is a
a cyclic permutation $\bb'$ of $\bb$ with $g(\bb')=1$.  
Thus, we expect that $1/f(\bb')$ will be $\gg 1/k$ on average over $\bb'$
and that $1/f(\bb) \gg 1/k$ on average over $\bb\in \BB$.  This
is essentially what we prove next; see \eqref{S0} below.

\begin{lem}\label{cycles}
For positive real numbers $x_1,\ldots,x_r$ with product $X$, let
$x_{r+i}=x_i$ for $i\ge 1$.  Then
$$
\sum_{j=0}^{r-1} \( \sum_{h=1}^r x_{1+j}\cdots x_{h+j} \)^{-1} \in 
\left[ \frac{1}{\max(1,X)}, \frac{1}{\min(1,X)} \right].
$$
\end{lem}

\begin{proof}
Put $y_0=1$ and $y_j=x_1\cdots x_j$ for $j\ge 1$.  The sum in question is
$$
\sum_{j=0}^{r-1} \( \sum_{h=1}^r \frac{y_{h+j}}{y_j} \)^{-1} = 
\sum_{j=0}^{r-1} \frac{y_j}{y_{1+j}+\cdots + y_{r+j}}.
$$
Since $y_r=X$,
\begin{align*}
y_{1+j}+\cdots + y_{r+j} &= X(y_0+\cdots+y_j)+y_{1+j}+\cdots+y_{r-1} \\
&\in [\min(1,X)(y_0+\cdots+y_{r-1}),\max(1,X)(y_0+\cdots+y_{r-1})]. \qedhere
\end{align*}
\end{proof}

We have
\be\label{low2}
\sum_{\bb\in\BB^*}  \frac{1}{b_M!\cdots b_J! f(\bb)} 
\ge S_0 - \sum_{M \le j < k/M} S_1(j) - \sum_{1\le j <k/M} S_2(j),
\ee
where
\begin{align*}
S_0 &= \sum_{\bb\in\BB} \frac{1}{b_M!\cdots b_J! f(\bb)}, \\
S_1(j) &= \sum_{\substack{\bb\in\BB \\ b_j>Mj}} 
  \frac{1}{b_M!\cdots b_J! f(\bb)}, \\
S_2(j) &= \sum_{\substack{\bb\in\BB \\ b_{J+1-j}>Mj}} 
  \frac{1}{b_M!\cdots b_J! f(\bb)}.
\end{align*}
Let $x_i=2^{-1+b_{M-1+i}}$ for $1\le i\le k$.  Then
$x_1\cdots x_k=1$ and
$$
f(\bb)=x_1 + x_1 x_2 + \cdots + x_1x_2\cdots x_k.
$$
By Lemma \ref{cycles} and the multinomial theorem,
\be\label{S0}
S_0 = \sum_{\bb\in\BB} \frac{1}{b_M!\cdots b_J!} \; \frac{1}{k}
\sum_{j=0}^{k-1} \( \sum_{h=1}^k x_{1+j}\cdots x_{h+j} \)^{-1} 
= \frac{k^{k-1}}{k!}.
\ee
To bound $S_1(j)$, apply Lemma \ref{cycles} with $x_i=2^{b_{j+i}-1}$ for
$1\le i\le J-j$ and note that 
$$
X=x_1\cdots x_{J-j}=2^{j+1-M-b_M-\cdots-b_j}<1.
$$
Write $\bb'=(b_M,\ldots,b_{j-1},b_{j+1},\ldots,b_J)$,
whose sum of components is $k-b_j$.  Ignoring the terms with $h\le j$ in
\eqref{fdef}, using Lemma \ref{cycles} and the multinomial theorem, we find
\begin{align*}
S_1(j) &\le \sum_{b_j>Mj} \frac{1}{b_j!} 
\sum_{\bb'} \frac{1}{\prod_{i\ne j} b_i!} 
\frac{1}{2^{M-1-j+b_M+\cdot+b_j}} \frac{1}{J-j}
\sum_{i=0}^{J-j-1} \( \sum_{h=1}^{J-j} x_{1+i}\cdots x_{h+i} \)^{-1} \\
&\le \frac{1}{J-j} \sum_{b_j>Mj}
\frac{(k-1)^{k-b_j}}{b_j! (k-b_j)!}
\le \frac{2 k^{k-1}}{k!} \sum_{b_j>Mj} \frac{1}{b_j!}
\le \frac{k^{k-1}}{k!} \frac{2}{(Mj)!}.
\end{align*}
Hence, if $M\ge 2$ then
\be\label{S1}
\sum_{M \le j < k/M} S_1(j) \le \frac{k^{k-1}}{10k!}.
\ee
The estimation of $S_2(j)$ is similar.  Let $x_i=2^{-1+b_{M+i-1}}$
for $1\le i\le J-M+1-j$, so that 
$$
X=x_1\cdots x_{J-M+1-j} = 2^{j-b_{J-j+1}-\cdots-b_J} \le 1.
$$
Put $b=b_{J-j+1}$ and let $\bb'=(b_M,\ldots,b_{J-j},b_{J-j+2},\ldots,b_J)$,
whose sum of components is $k-b$.  Then, ignoring the terms with
$h>J-j$ in \eqref{fdef}, we have
\begin{align*}
S_2(j) &\le \sum_{b>Mj} \frac{1}{b!} \sum_{\bb'} \frac{1}
{\prod_{i\ne J-j+1} b_i!} \frac{2^{b-j+b_{J-j+2}+\cdots+b_J}}{J-M+1-j} \\
&= \frac{2^{-j}}{J-M+1-j} 
  \sum_{b>Mj} \frac{2^b}{b!} \; \frac{(k+j-2)^{k-b}}{(k-b)!} \\
&\le \frac{2^{1-j}}{k\cdot k!} (k+j)^k \sum_{b>Mj} \frac{2^b}{b!} \\
&\le \frac{k^{k-1}}{k!} 2^{1-j} e^j \frac{2^{Mj}}{(Mj)!}.  
\end{align*}
If $M$ is large enough, then
\be\label{S2}
\sum_{j\ge 1} S_2(j) \le \frac{k^{k-1}}{10k!}.
\ee
By \eqref{low2}, \eqref{S0}, \eqref{S1} and \eqref{S2},
$$
\sum_{\bb\in\BB^*}  \frac{1}{b_M!\cdots b_J! f(\bb)} 
\ge \frac{k^{k-1}}{2k!}.
$$
The lower bound in \eqref{Hxy2y} for large $y$ now follows from \eqref{low1} 
and Stirling's formula.  If $y\le y_0$ for some fixed constant $y_0$,
the lower bound in \eqref{Hxy2y} follows from $H(x,y,2y)\gg x$.

%
%
%
%
%
\section{Upper bound, part I}\label{sec:upper1}
%
%
%
%
%

In this section, we prove the upper bound implicit in \eqref{Hxy2y},
except for the estimation of some integrals which will be dealt with in
section \ref{sec:upper2}.  As with the lower bound argument, we begin by
bounding $H(x,y,2y)$ in terms of a sum involving $L(a)$.  Using a relatively
simple upper bound for $L(a)$ proved in Lemma \ref{L} below, the sums
involving $L(a)$ are bounded in terms of particular multivariate integrals.
The estimates for these integrals in section \ref{sec:upper2} allow us then to
complete the proof.

\begin{lem}\label{L}
We have
\begin{enumerate}
\item  $L(a) \le \min(\tau(a)\log 2, \log 2+\log a)$;
\item  If $(a,b)=1$, then $L(ab) \le \tau(b) L(a)$;
\item  If $p_1 < \cdots < p_k$, then
$$
L(p_1\cdots p_k) \le\min_{0 \le j\le k} 2^{k-j} (\log(p_1\cdots p_j)+\log 2).
$$
\end{enumerate}
\end{lem}

\begin{proof}  Part (i) is immediate, since $\LL(a)$ is the union of
$\tau(a)$ intervals of length $\log 2$, all contained in $[-\log 2,\log a)$.
Part (ii) follows from
\benn
\LL(ab) = \bigcup_{d|b} \{ u + \log d : u\in \LL(a)\}.
\eenn
Combining parts (i) and (ii) with $a=p_1\cdots p_j$ and
$b=p_{j+1}\cdots p_k$ yields (iii).
\end{proof}

\begin{lem}\label{HL}
If $3\le y\le \sqrt{x}$, then 
\[
H(x,y,2y) \ll x \max_{\sqrt{y} \le t \le x}  \sum_{\substack{P^+(a) \le t \\ \mu^2(a)=1}}  
\frac{L(a)}{a \log^2(t/a+P^+(a))}.
\]
\end{lem}

\begin{proof}  
First, we relate $H(x,y,2y)$ to $H^*(x,y,z)$, the number of
\emph{squarefree} integers $n\le x$ with $\tau(n,y,z)\ge 1$.  Write
$n=n'n''$, where $n'$ is squarefree, $n''$ is squarefull and $(n',n'')=1$.
The number of $n\le x$ with $n''>(\log y)^4$ is
$$
\le x \sum_{n''>(\log y)^4} \frac{1}{n''} \ll \frac{x}{(\log y)^2}.
$$
Assume now that $n'' \le (\log y)^4$.  For some $f|n''$, $n'$ has a divisor in
$(y/f,2y/f]$, hence
\be\label{HH*}
H(x,y,2y)\le \sum_{n''\le (\log y)^4} \sum_{f|n''} H^*\( \tfrac{x}{n''},
\tfrac{y}{f}, \tfrac{2y}{f} \) + O\pfrac{x}{(\log y)^2}.
\ee
Next, we show that for $3\le y_1\le x_1^{3/5}$,
\be\label{e2}
H^*(x_1,y_1,2y_1)-H^*(\tfrac12 x_1,y_1,2y_1) \ll x_1 
\bigl( S(2y_1)+S(x_1/y_1) \bigr),
\ee
where
\[
S(t) = \sum_{\substack{P^+(a) \le t \\ \mu^2(a)=1}}  
\frac{L(a)}{a \log^2(t/a+P^+(a))}.
\]
Let $\AA$ be the set of squarefree integers $n\in (\frac12 x_1,x_1]$ 
with a divisor in $(y_1,2y_1]$.
Put $z_1=2y_1$, $y_2=\frac{x_1}{4y_1}$, $z_2=\frac{x_1}{y_1}$.
If $n\in \AA$, then $n=m_1m_2$ with $y_i < m_i \le z_i$ ($i=1,2$).
For some $j\in \{1,2\}$ we have $p=P^+(m_j) < P^+(m_{3-j})$.
Write $n=abp$, where $P^+(a) < p < P^-(b)$ and $b>p$.  
Since $\tau(ap,y_j,z_j)\ge 1$, we have $p \ge y_j/a$.
By \eqref{sieve},
given $a$ and $p$, the number of possible $b$ is
$$
\ll \frac{x_1}{ap\log p} \le \frac{x_1}{ap\log \max(P^+(a),y_j/a)},
$$
Since $a$ has a divisor in $(y_j/p,z_j/p]$,
we have $\log(y_j/p)\in \LL(a)$ or $\log(2y_j/p)\in \LL(a)$.  Since
$\LL(a)$ is the disjoint union of intervals of length $\ge \log 2$
with total measure $L(a)$, by repeated use of \eqref{sum1p} we obtain
$$
\sum_{\substack{\log (cy_j/p)\in \LL(a) \\ p\ge P^+(a)}} \frac{1}{p} \ll
\frac{L(a)}{\log \max(P^+(a),y_j/a)} \qquad (c=1,2),
$$
and \eqref{e2} follows.

Write $x_2=x/n''$, $y_1=y/f$.
Each $n\in (x_2/\log^2 y_1,x_2]$ lies in an interval $(2^{-r+1}x_2,2^{-r}x_2]$ 
 for some integer $0\le r\le 5\log\log y_1$.  
Applying \eqref{e2} with $x_1=2^{-r}x_2$ for each $r$ gives
$$
H^*(x_2,y_1,2y_1)\ll \frac{x_2}{\log^2 y_1} + \sum_r 2^{-r} x_2 
\bigl( S(2y_1)+S(2^{-r}x_2/y_1) \bigr) 
\ll x_2 \max_{\sqrt{y_1} \le t\le x_2} S(t).
$$
Here we used the fact that $S(t) \ge \frac{L(1)}{\log^2 t} = \frac{\log
  2}{\log^2 t}$. 
Finally, $\sum_{n''} \tau(n'')/n'' = O(1)$ and 
the lemma follows from \eqref{HH*}.
\end{proof}

The next lemma is due to Kouloulopoulos \cite[Lemma 2.2]{DK}.  We give
 a much shorter proof.

\begin{lem}\label{lem22}
Suppose $f$ is an arithmetic function satisfying $f(pm)\le Cf(m)$ for
all primes $p$ and all $m\in\NN$ coprime to $p$.  Let
$\PP(x) = \{n\in \NN : \mu^2(n)=1, P^+(n)\le x \}$.
For any real $h\ge 0$,
\[
\sum_{a\in\PP(x)}
\frac{f(a)}{a\log^h(P^+(a)+x/a)}
\ll_{C,h} \frac{1}{(\log x)^h} \sum_{a\in\PP(x)} \frac{f(a)}{a}.
\]
\end{lem}

\begin{proof}
Let $\PP_1=\{a\in\PP(x):a>x^{1/2},P^+(a)\le x^{1/4}\}$.  Then clearly
\[
\sum_{a\in\PP(x)}
\frac{f(a)}{a\log^h(P^+(a)+x/a)} \le
\sum_{a\in\PP_1} \frac{f(a)}{a\log^h P^+(a)} +
\frac{4^h}{(\log x)^h}  \sum_{a\in\PP(x)} \frac{f(a)}{a}.
\]
For $a\in\PP_1$, let $p=P^+(a)$ and $a=pb$, so $b>x^{1/4}$.  
Let $k=\lfloor h+2 \rfloor$. Since $f(pb)\le Cf(b)$,
\begin{align*}
\sum_{a\in\PP_1} \frac{f(a)}{a\log^h P^+(a)} &\le
C \sum_{p\le x^{1/4}} \frac{1}{p\log^h p} \sum_{\substack{b\in\PP(p)
      \\ b > x^{1/4}}} \frac{f(b)}{b} \\
&\le C \sum_{p\le x^{1/4}} \frac{1}{p\log^h p} \,
  \frac{4^k}{\log^{k} x} \sum_{b\in\PP(p)} \frac{f(b)\log^{k} b}{b}.
\end{align*}
Next,
\begin{align*}
\sum_{b\in\PP(p)} \frac{f(b)\log^{k} b}{b} &= \sum_{b\in\PP(p)}
\frac{f(b)}{b} \sum_{p_1|b, \ldots,p_{k}|b} (\log p_1)\cdots (\log p_k) \\
&\le C^k \sum_{p_1,\ldots,p_k \le p} \frac{(\log p_1)\cdots (\log
  p_k)}{[p_1,\ldots,p_k]} \sum_{t\in \PP(p)} \frac{f(t)}{t},
\end{align*}
where we have written $b=[p_1,\ldots,p_k]t$. 
Write $p_1\cdots p_k = q_1^{e_1}\cdots
q_m^{e_m}=r$, where $q_1 < \cdots< q_m$ are prime.  With $r$ fixed,
there are $O_k(1)$ choices for $p_1,\ldots,p_k$.  Also, there are
$O_k(1)$ choices for $e_1,\ldots,e_m$ for each choice of $m$.  Hence, the sum on
$p_1,\ldots,p_k$ is
\[
\ll_h \sum_{m=1}^k \;\; \sum_{e_1+\cdots+e_m=k} \;\; \prod_{j=1}^m \Biggl(
\sum_{q<p} \frac{(\log q)^{e_j}}{q} \Biggr) \ll_h (\log p)^k
\]
by repeated application of Mertens' estimate.  Extending the range of 
$t$ to $t\in \PP(x)$, we get
\[
\sum_{a\in\PP_1} \frac{f(a)}{a\log^h(P^+(a)+x/a)} \ll_{C,h}
\sum_{p\le x^{1/4}} \frac{1}{p\log^h p}\, \frac{(\log p)^k}{(\log
  x)^k} \sum_{t\in \PP(x)} \frac{f(t)}{t}.
\]
A final application of Mertens' estimate concludes the proof.
\end{proof}

By Lemma \ref{L}, the hypotheses of Lemma \ref{lem22} are satisfied with
$f(a)=L(a)$ and $C=2$.
Combining Lemma \ref{HL} with \ref{lem22} produces
an upper bound of the same shape as the lower bound in Lemma \ref{HL_lower}.

\begin{lem}\label{HL2}
Uniformly for $3\le y\le \sqrt{x}$, we have
\[
H(x,y,2y) \ll x \max_{\sqrt{y} \le t \le x}  \frac{1}{\log^2 t}
\sum_{\substack{P^+(a) \le t \\ \mu^2(a)=1}}  \frac{L(a)}{a}.
\]
\end{lem}

We cut up the sum in Lemma \ref{HL2} according to $\omega(a)$.  Let
$$
T_k(P) = \sum_{\substack{P^+(a)\le P, \, \mu^2(a)=1 \\ \omega(a)=k}}
\frac{L(a)}{a}.
$$

%
%

We next bound $T_k(P)$ in terms of a mutivariate integral.  Since
$\sum_{p\le z} 1/p = \log\log z+O(1)$, by partial summation we expect 
for ``nice'' functions $f$ that
$$
\sum_{p_1 < \cdots < p_k \le P} \frac{f\( \frac{\log\log p_1}{\log\log P},
 \cdots,\frac{\log\log p_k}{\log\log P}\)}{p_1\cdots p_k} 
\approx (\log\log P)^k 
\idotsint\limits_{0\le \xi_1 \le \cdots\le \xi_k\le 1} f(\bx)\, d\bx.
$$

\begin{lem}\label{Tnint}
Suppose $P$ is large, $v=\fl{\frac{\log\log P}{\log 2}}$
and $1\le k\le 10v$.  Then
\[
T_k(P) \ll (2\log\log P)^k U_k(v), \quad
U_k(v)= \idotsint\limits_{0\le \xi_1 \le \cdots\le \xi_k\le 1} 
\min_{0\le j\le k} 2^{-j} (2^{v\xi_1}+\cdots+2^{v\xi_j}+1)\, d\bx.
\]
\end{lem}

\begin{proof}
Recall the definition of $\lambda_i, D_i$ from \S 2.
Consider $a=p_1 \cdots p_k$, $p_1 < \cdots < p_k \le P$ 
and define $j_i$ by 
$p_i\in D_{j_i}$  ($1\le i\le k$).  Put $l_i = \frac{\log\log p_i}{\log 2}$.
By Lemma \ref{L} (iii) and \eqref{lam},
$$
L(a) \le 2^k \min_{0\le g\le k} 2^{-g} (2^{l_1}+\cdots +2^{l_g}+1)
\le 2^{k+K} F(\bj),
$$
where
$$
F(\bj) = \min_{0\le g\le k} 2^{-g} (2^{j_1}+\cdots+ 2^{j_g}+1).
$$
Let $J$ denote the set of vectors $\bj$ satisfying
$0\le j_1 \le \cdots \le j_k \le v+K+1.$  Then
$$
T_k(P) \le 2^{k+K} \sum_{\bj \in J} F(\bj) 
\sum_{\substack{p_1< \cdots < p_k \\ p_i\in D_{j_i} \; (1\le i\le k)}} 
\frac{1}{p_1\cdots p_k}.
$$
Let $b_j$ be the number of primes $p_i$ in $D_{j}$ for $0\le j\le v+K-1$.
Using the hypothesis that $k\le 10v$,
the sum over $p_1,\cdots,p_k$ above is at most
\begin{align*}
\prod_{j=0}^{v+K+1} \frac{1}{b_j!} \biggl( \sum_{p\in E_j}
\frac{1}{p} \biggr) ^{b_j}
&\le \frac{(\log 2)^k}{b_0! \cdots b_{v+K+1}!} \\
&= ((v+K)\log 2)^k
\int_{R(\bj)} 1\, d\bx 
\le e^{10K} (v\log 2)^k \int_{R(\bj)} 1\, d\bx,
\end{align*}
where
$$
R(\bj) = \{ 0\le \xi_1\le \cdots \le \xi_k\le 1: 
j_i \le (v+K+2) \xi_i \le j_i+1 \,\; \forall i\}.
$$
In $R(\bj)$, there are $b_s$ numbers $\xi_j$ satisfying
$s\le (v+K+2)\xi_i \le s+1$ for each $s$, and $\text{Vol} \{0\le x_1 \le \cdots
\le x_b\le 1\} = 1/b!$.
Since $2^{j_i} \le 2^{(v+K+2)\xi_i} \le 2^{K+2} 2^{v\xi_i}$ for $\bx\in R(\bj)$, we
have 
$$
F(\bj) \le 2^{K+2} \min_{0\le g\le k} 2^{-g}(2^{v\xi_1}+\cdots+2^{v\xi_g}+1).
$$ 
Hence
$$
\sum_{\bj \in J} F(\bj) \int_{R(\bj)} 1 d\bx \le 2^{K} U_k(v)
$$
and the lemma follows.
\end{proof}

Estimating $U_k(v)$ is the most complex part of the argument.
The next lemma will be proved in section \ref{sec:upper2}.

\begin{lem}\label{Unlem}
Suppose $k,v$ are integers with $0\le k\le 10v$.  Then
$$
U_k(v) \ll \frac{1 + |v-k|^2} {(k+1)! (2^{k-v}+1)}.
$$
\end{lem}

Notice that the bound in Lemma \ref{Unlem} undergoes a change of behavior
at $k = v$.  

\begin{proof}[Proof of \eqref{Hxy2y}, upper bound]
Let $v=\fl{\frac{\log\log P}{\log 2}}$.  By Lemmas \ref{Tnint}
and \ref{Unlem},
\[
\sum_{v\le k\le 10v} T_k(P) \ll  \sum_{v\le k\le 10v}
\frac{(k-v)^2+1}{2^{k-v}} \; \frac{(2\log\log P)^{k}}{(k+1)!} 
\ll \frac{(2\log\log P)^v}{(v+1)!}
\]
and
\[
\sum_{1\le k\le v} T_k(P) \ll \sum_{1\le k\le v}
\frac{((v-k)^2+1)(2\log\log P)^k}{(k+1)!} \ll \frac{(2\log\log
  P)^v}{(v+1)!}. 
\]
By Lemma \ref{L} (i),
\begin{align*}
\sum_{k\ge 10v} T_k(P) &\le \sum_{k\ge 10v} \sum_{\substack{P^+(a) \le P \\
    \mu^2(a)=1, \omega(a)=k}} \frac{2^k \log 2}{a}
\le \sum_{k\ge 10v} \frac{2^k}{k!} \biggl( \sum_{p\le P}
    \frac{1}{p} \biggr)^k \\
&\le \frac{(2\log\log P+O(1))^{10v}}{(10v)!}
 \ll \frac{(2\log\log P)^v}{(v+1)!}.
\end{align*}
Finally, $T_0(P)=L(1)=\log 2$.
Recalling the definition of $v$ and combining the above bounds on
$T_k(P)$ with Stirling's formula and Lemma \ref{HL2} completes the proof.
\end{proof}

%
%
%
\section{Upper bound, part II}\label{sec:upper2}
%
%
%

The goal of this section is to prove Lemma \ref{Unlem}, and thus complete the
proof of the upper bound in \eqref{Hxy2y}.

Let $Y_1, \ldots, Y_n$ be independent, uniformly distributed random
variables in $[0,1]$.  Let $\xi_1$ be the smallest of the numbers
$Y_i$, let $\xi_2$ be the next smallest, etc., so that
$0 \le \xi_1 \le \cdots \le \xi_n \le 1.$
The numbers $\xi_i$ are the \emph{order statistics} for 
$Y_1, \ldots, Y_n$.  Then $k! U_k(v)$ is the expectation of the random
variable 
$$
X=\min_{0\le j\le k} 2^{-j}(2^{v\xi_1}+\cdots+2^{v\xi_j}+1).
$$
Heuristically, we expect that 
\be\label{EX}
\EE X\ll \EE \min_{1\le j\le k} 2^{-j+v\xi_j}, 
\ee
so we need to understand the distribution of
$\min_{1\le j\le k} v\xi_j-j$.  Let $Q_k(u,v)$ be
the probability that $\xi_i \ge \frac{i-u}{v}$ for every $i$.
In the special case $v=k$, Smirnov in 1939 showed that
$$
Q_k(x\sqrt{k},k) \sim 1 - e^{-2x^2}
$$ 
for each fixed $x$.  The corresponding probability estimate for two-sided
bounds on the $\xi_i$ was established by Kolmogorov in 1933 and together these
limit theorems are the basis of the \emph{Kolmogorov-Smirnov goodness-of-fit
  statistical tests}.

In the next lemma, we prove new, uniform estimates for $Q_k(u,v)$.
The remainder of the section is essentially devoted to proving \eqref{EX}.
The details are complicated, but the basic idea is that if $2^{-j}(2^{v\xi_1}+
\cdots + 2^{v\xi_j})$ is much larger than $2^{v\xi_j-j}$, then for some large
$l$, the numbers $\xi_{j-l},\ldots,\xi_j$ are all very close to one another.
As shown below in Lemmas \ref{vol1} and \ref{UUlem}, this is quite rare.

\begin{lem}\label{Q}  Let $w=u+v-k$.  Uniformly in $u \ge 0$ and $w \ge 0$, 
we have
$$
Q_k(u,v) \ll \frac{(u+1)(w+1)^2}{k}.
$$
\end{lem}
\begin{proof}
Without loss of generality, suppose $k\ge 100$, $u\le k/10$ and 
$w \le \sqrt{k}$.  If $\min_{1\le i\le k} (\xi_i
- \frac{i-u}{v})<0$, let $l$ be the smallest index with $\xi_l <
\frac{l-u}{v}$ and write $\xi_l=\frac{l-u-\lam}{v}$, so that $0\le \lam
\le 1$.  Let
$$
R_l(\lam) = \Vol\left\{ 0\le \xi_1 \le \cdots \le \xi_{l-1} \le
\frac{l-u-\lam}{v} : \xi_i \ge \frac{i-u}{v} \, (1\le i\le l-1) \right\}.
$$
Then we have
\begin{align*}
Q_k(u,v) &= 1 - \frac{k!}{v} \int_0^1 \sum_{u+\lam \le l \le k}
R_l(\lam) \Vol \left\{ \frac{l-u-\lam}{v} \le \xi_{l+1} \le \cdots \le
\xi_k \le 1 \right\} \, d\lam \\
&= 1 - \frac{k!}{v} \int_0^1 \sum_{u+\lam \le l \le k}
\frac{R_l(\lam)}{(k-l)!} \pfrac{k+w+\lam-l}{v}^{k-l}\, d\lam.
\end{align*}

Now suppose that $\xi_k \le 1 - \frac{2w+2}{v} = \frac{k-u-w-2}{v}$.
Then $\min_{1\le i\le k} (\xi_i-\frac{i-u}{v})<0$.  Defining $l$ and $\lam$
as before, we have
\begin{align*}
\(1 - \frac{2w+2}{v}\)^k &= k! \Vol \biggl\{ 
0\le \xi_1\le \cdots \le \xi_k \le 1 -
\frac{2w+2}{v} \biggr\} \\
&= \frac{k!}{v} \int_0^1 \sum_{u+\lam \le l \le k-w-2+\lam} 
\frac{R_l(\lam)}{(k-l)!} \pfrac{k-l-w-2+\lam}{v}^{k-l}\, d\lam.
\end{align*}
Thus, for any $A>0$, we have
\begin{multline*}
Q_k(u,v) = 1 - A \( 1 - \frac{2w+2}{v} \)^k - \frac{k!}{v} \int_0^1
\sum_{k-w-2+\lam <l \le k} \frac{R_{l}(\lam)}{(k-l)!} 
\pfrac{k+w+\lam-l}{v}^{k-l}\, d\lam \\
+ \frac{k!}{v} \int_0^1 \sum_{u+\lam \le l \le k-w-2+\lam}
  \frac{R_{l}(\lam)}{(k-l)! v^{k-l}} \left[
    A(k-l-w-2+\lam)^{k-l} - (k-l+w+\lam)^{k-l} \right]\, d\lam. 
\end{multline*}
Noting that $2-\lam \ge \lam$, we have
\begin{align*}
\pfrac{k-l-w-2+\lam}{k-l+w+\lam}^{k-l} &= \(1 - \frac{w+2-\lam}{k-l}\)^{k-l}
\(1 + \frac{w+\lam}{k-l}\)^{-(k-l)} \\
&= \exp \left\{ -(2w+2) +
\sum_{j=2}^\infty \frac{-(w+2-\lam)^j + (-1)^j(w+\lam)^j}{j(k-l)^{j-1}}
\right\} \\
&\le e^{-(2w+2)}.
\end{align*}
Thus, taking $A=e^{2w+2}$, we conclude that
\begin{align*}
Q_k(u,v) &\le 1 - e^{2w+2} \(1 - \frac{2w+2}{v}\)^k \\
&= 1 - \exp \left\{ \frac{2w+2}{v} \( v - k + O(w)\) 
  \right\} \\
&= 1 - \exp \left\{ \frac{-2uw+O(u+w^2+1)}{v} \right\} \\
&\le \frac{2uw+O(u+w^2+1)}{v} \ll \frac{(u+1)(w+1)^2}{k}.\qedhere 
\end{align*}
\end{proof}

\begin{lem}\label{combsum} 
If $t\ge 2$, $b\ge 0$ and $t+a+b>0$, then
$$
\sum_{\substack{1\le j\le t-1 \\ j+a>0}} \binom{t}{j} (a+j)^{j-1}
(b+t-j)^{t-j-1} \le e^{4} (t+a+b)^{t-1}.
$$
\end{lem}

\begin{proof} 
Let $C_t(a,b)$ denote the sum in the lemma.  
We may assume that $a>1-t$, otherwise $C_t(a,b)=0$.  
The associated ``complete'' sum is evaluated
exactly using one of Abel's
identities (\cite{Riordan}, p.20, equation (20))
\be\label{Abel}
\sum_{j=0}^t \binom{t}{j} (a+j)^{j-1} (b+t-j)^{t-j-1} = \( \frac{1}{a} +
\frac{1}{b} \) (t+a+b)^{t-1} \qquad (ab\ne 0).
\ee
If $a\ge -1$, put $A=\max(1,a)$ and $B=\max(1,b)$.  By
\eqref{Abel},
\be\label{Ctab}
\begin{split}
C_t(a,b) &\le C_t(A,B) \le \( \frac{1}{A} + \frac{1}{B} \) (t+A+B)^{t-1} \\
&\le 2 (t+a+b+3)^{t-1} \\
&\le 2 e^{\frac{3(t-1)}{t+a+b}} (t+a+b)^{t-1} < e^4 (t+a+b)^{t-1}. 
\end{split}
\ee
Next assume $a< -1$.  Since $(1+c/x)^x$ is an increasing function 
for $x>1$, we have
\benn
(a+j)^{j-1}=(j-1)^{j-1} \( 1 + \frac{a+1}{j-1} \)^{j-1} \le (j-1)^{j-1}
\( 1 + \frac{a+1}{t-1} \)^{t-1}.
\eenn
Thus, by \eqref{Ctab},
\begin{align*}
C_t(a,b) &\le \pfrac{t+a}{t-1}^{t-1}  C_t(-1,b) \\
&\le e^4 \pfrac{(t+a)(t+b-1)}{t-1}^{t-1} = e^4 \(t+a+b + 
  \frac{(a+1)b}{t-1} \)^{t-1} \\
&\le e^{4} (t+a+b)^{t-1}. \qedhere
\end{align*}
\end{proof}

For brevity, write
$$
S_k(u,v) = \{ \bx : 0 \le \xi_1 \le \cdots \le \xi_k \le 1 : \xi_i \ge
\frac{i-u}{v} \, (1\le i\le k) \}, 
$$
so that $Q_k(u,v) = k! \Vol S_k(u,v)$.

\begin{lem}\label{vol1}
Suppose $g,k,s,u,v \in \ZZ$ satisfy
$$
1\le g\le k-1, \; s\ge 0, \; v\ge k/10, \; u\ge 0, \; u+v\ge k+1.
$$
Let $R$ be the subset of $\bx \in S_k(u,v)$ where, for some $l\ge g+1$,
we have
\be\label{lcond1}
\frac{l-u}{v} \le \xi_l \le \frac{l-u+1}{v}, \qquad \xi_{l-g} \ge 
\frac{l-u-s}{v}.
\ee
Then
$$
\Vol(R) \ll \frac{g^2(10(s+1))^g}{g!} \, \frac{(u+1)(u+v-k)^2}{(k+1)!}.
$$
\end{lem}

\begin{proof}  Fix $l$ satisfying $\max(u,g+1) \le l \le k$.  Let $R_l$
be the subset of  $\bx \in S_k(u,v)$ satisfying \eqref{lcond1} for
this particular $l$. 
We have $\Vol(R_l) \le V_1 V_2 V_3 V_4$, where, by Lemma \ref{Q},
\begin{align*}
V_1 &= \Vol\{ 0\le \xi_1 \le \cdots \le \xi_{l-g-1} \le \tfrac{l-u+1}{v}:
  \xi_i \ge \tfrac{i-u}{v}\, \forall i \} \\
&= \pfrac{l-u+1}{v}^{l-g-1} 
  \Vol\{ 0\le \theta_1 \le \cdots \le \theta_{l-g-1}\le 1 : \theta_i \ge
  \tfrac{i-u}{l-u+1}\; \forall i \}  \\
&= \pfrac{l-u+1}{v}^{l-g-1} \frac{Q_{l-g-1}(u,l-u+1)}{(l-g-1)!} \\
&\ll \pfrac{l-u+1}{v}^{l-g-1} \frac{(u+1) g^2}{(l-g)!},
\end{align*}
\begin{align*}
V_2 &= \Vol \{ \tfrac{l-u-s}{v} \le \xi_{l-g} \le \cdots \le \xi_{l-1} \le
  \tfrac{l-u+1}{v} \} = \frac{1}{g!} \pfrac{s+1}{v}^g, \\
V_3 &= \Vol \{ \tfrac{l-u}{v} \le \xi_l \le \tfrac{l-u+1}{v} \} = 
  \frac{1}{v}, \\
V_4 &= \Vol\{ \xi_{l+1} \le \cdots \le \xi_k \le 1 : \xi_i \ge \tfrac{i-u}{v}
  \, \forall i \} \\
&= \frac{1}{(k-l)!} \pfrac{u+v-l}{v}^{k-l} Q_{k-l}(0,u+v-l) \\
&\ll \pfrac{u+v-l}{v}^{k-l} \frac{(u+v-k)^2}{(k-l+1)!}.
\end{align*}
Thus
$$
\Vol(R) \ll \frac{(s+1)^g (u+1) g^2 (u+v-k)^2}{g! v^k (k+1-g)!}
\sum_l \binom{k+1-g}{l-g} (l-u+1)^{l-g-1} (u+v-l)^{k-l}.
$$
By Lemma \ref{combsum} (with $t=k+1-g$, $a=g+1-u$, $b=u+v-k-1$), 
the sum on $l$ is 
$$
\le e^4 (v+1)^{k-g} 
\ll v^{k-g} = \frac{v^k}{k^g} \,  \pfrac{k}{v}^g \le v^k 10^g
\frac{(k-g+1)!}{k\cdot k!}
$$
and the lemma follows.
\end{proof}

To bound $U_k(v)$, we will bound the volume of the set
$$
\TT(k,v,\g) = \{ \bx\in \RR^k : 0 \le \xi_1 \le \cdots \le \xi_k \le 1, 
  2^{v\xi_1} + \cdots + 2^{v\xi_j} \ge 2^{j-\g}\; (1\le j\le k) \}.
$$

\begin{lem}\label{UUlem}
Suppose $k,v,\g$ are integers with $1\le k\le 10v$ and $\g\ge 0$.  
Set $b=k-v$.  Then
$$
\Vol(\TT(k,v,\g)) \ll \frac{Y}{2^{2^{b-\g}} (k+1)!}, \qquad
Y = \begin{cases} b & \text{ if } b\ge \g+5 \\ (\g+5-b)^2(\g+1) & 
\text{ if } b\le \g+4  \end{cases}.
$$
\end{lem}

\begin{proof}  
Let $r=\max(5,b-\g)$ and $\bx \in \TT(k,v,\g)$.  Then either
\be\label{U-A1}
\xi_j > \tfrac{j-\g-r}{v} \quad (1\le j\le k)
\ee
or
\be\label{U-A2}
\min_{1\le j\le k} ( \xi_j - \tfrac{j-\g}{v} ) = \xi_l - \tfrac{l-\g}{v} \in 
[\tfrac{-h}{v}, \tfrac{1-h}{v} ] \; \text{ for some integers } h\ge
r+1, 1\le l \le k.
\ee
Let $V_1$ be the volume of
$\bx \in \TT(k,v,\g)$ satisfying \eqref{U-A1}.
If $b\ge \g+5$, \eqref{U-A1} is not possible, so $b\le \g+4$ and $r=5$.
By Theorem \ref{Q},
$$
V_1 \le \frac{Q_k(\g+5,v)}{k!} \ll \frac{(\g+6)(\g+6-b)^2}{(k+1)!}
\ll \frac{Y}{2^{2^{b-\g}} (k+1)!}.
$$

If \eqref{U-A2} holds, then there is an integer $m$ satisfying
\be\label{U-A3}
m\ge h-3, \; \xi_{l-2^m} \ge \tfrac{l-\g-2m}{v}.
\ee
To see \eqref{U-A3}, suppose such an $m$ does not exist.  Then
$$
2^{v\xi_1}+\cdots+2^{v\xi_l} \le 
 2^{h-3} 2^{l-\g-h+1} + \sum_{m\ge h-3} 2^m 2^{l-\g-2m}\le 2^{l-\g},
$$
a contradiction.  
Let $V_2$  be the volume of
$\bx \in \TT(k,v,\g)$ satisfying \eqref{U-A2}.
Fix $h$ and $m$ satisfying \eqref{U-A3}
and use Lemma \ref{vol1} with $u=\g+h$,
$g=2^m$, $s=2m$.  The volume of such $\bx$ is
$$
\ll \frac{(\g+h+1)(\g+h-b)^2}{(k+1)!}\, \frac{(20m+10)^{2^m}2^{2m}}{(2^m)!}
\ll \frac{(\g+h+1)(\g+h-b)^2}{2^{2^{m+3}} (k+1)!}.
$$
The sum of $2^{-2^{m+3}}$ over $m\ge h-3$ is $\ll 2^{-2^{h}}$.
Summing over $h\ge r+1$ gives
$$
V_2 \ll \frac{(\g+r+2)(\g-b+r+2)^2}{2^{2^{r+1}} (k+1)!}
\ll \frac{Y}{2^{2^{b-\g}} (k+1)!}.
$$
\end{proof}

\medbreak
\emph{Proof of Lemma \ref{Unlem}.}
Assume $k\ge 1$, since the lemma is trivial when $k=0$.
Put $b=k-v$.
For integers $m\ge 0$,  consider $\bx \in R_k$ satisfying
$2^{-m} \le\min_{0\le j\le k} 2^{-j}\( 2^{v\xi_1}+\cdots+2^{v\xi_j}+1
\)  < 2^{1-m}.$  For $1\le j\le k$ we have
$$
2^{-j} \( 2^{v\xi_1}+\cdots+2^{v\xi_j} \) \ge
\max(2^{-j}, 2^{-m}-2^{-j}) \ge 2^{-m-1},
$$
so $\bx \in \TT(k,v,m+1)$.  Hence, by Lemma \ref{UUlem},
\begin{align*}
U_k(v)  &\le \sum_{m\ge 0} 2^{1-m} \Vol(\TT(k,v,m+1))
\ll \frac{1}{(k+1)!} \sum_{m\ge 0} \frac{2^{-m} Y_m}
  {2^{2^{b-m-1}}}, \\
Y_m &= \begin{cases} b & \text{ if } m\le b-6 \\ (m+6-b)^2(m+2) & 
\text{ if } m\ge b-5 \end{cases}.
\end{align*}
Next,
$$
\sum_{m\ge 0} \frac{2^{-m} Y_m}{2^{2^{b-m-1}}} = \!\!\!
\sum_{0\le m\le b-6} \frac{b}{2^m 2^{2^{b-m-1}}} + \!\! \sum_{m\ge 
\max(0,b-5)}\!\! \frac{(m+6-b)^2(m+2)}{2^m}.
$$
The proof is completed by noting that
if $b\ge 6$, each sum on the right side is $\ll b 2^{-b}$ and
if $b\le 5$, the first sum is empty and the second is $\ll (6-b)^2 \ll 1+b^2$.
\qed

%
%

\bibliographystyle{amsplain}
\bibliography{hxy2y}

\end{document}